\documentclass{amsart}
\usepackage{math}
\begin{document}
\title{The Growth of Grigorchuk's Torsion Group}
\author{Laurent Bartholdi}
\date\today
\email{Laurent.Bartholdi@math.unige.ch}
\address{\parbox{.4\linewidth}{Section de Math\'ematiques\\
    Universit\'e de Gen\`eve\\ CP 240, 1211 Gen\`eve 24\\
    Switzerland}}
\thanks{This work has been supported by the ``Swiss National Fund for Scientific Research''}
\begin{abstract}
  In 1980 Rostislav Grigorchuk constructed a group $G$ of intermediate
  growth, and later obtained the following estimates on its
  growth $\gamma$~\cite{grigorchuk:gdegree}:
  $$e^{\sqrt{n}}\precsim\gamma(n)\precsim e^{n^\beta},$$
  where $\beta=\log_{32}(31)\approx0.991$. Using elementary methods we
  improve the upper bound to
  $$\gamma(n)\precsim e^{n^\alpha},$$
  where $\eta\approx0.811$ is the real root of the polynomial
  $X^3+X^2+X-2$ and $\alpha=\log(2)/\log(2/\eta)\approx0.767$.
\end{abstract}
\maketitle

\section{Introduction}
The notion of growth for finitely generated groups was introduced in
the 1950's in the former \textsc{Ussr}~\cite{svarts:growth} and in the
1960's in the West~\cite{milnor:solvable}.  There are well-known
classes of groups of polynomial growth (abelian, and more generally
virtually nilpotent groups~\cite{gromov:nilpotent}) and of exponential
growth (non-virtually-nilpotent linear~\cite{tits:linear} or
non-elementary hyperbolic~\cite{ghys-h:gromov} groups).  However, the
first example of a group of intermediate growth was discovered later,
by Rostislav Grigorchuk;
see~\cite{grigorchuk:growth,grigorchuk:gdegree,grigorchuk:kyoto}. He
showed that the growth $\gamma$ of his group satisfies
$$e^{\sqrt{n}}\precsim\gamma(n)\precsim e^{n^\beta},$$
where $\beta=\log_{32}(31)\approx0.991$; see below for the precise
definition of growth.
The purpose of this note is to prove the following improvement:
\begin{thm}\label{thm:main}
  Let $\eta$ be the real root of the polynomial $X^3+X^2+X-2$, and set
  $\alpha=\log(2)/\log(2/\eta)\approx0.767$. Then the growth $\gamma$ of
  Grigorchuk's group satisfies
  $$e^{\sqrt{n}}\precsim\gamma(n)\precsim e^{n^\alpha}.$$
\end{thm}

\section{Growth of Groups}
Let $G$ be a group generated as a monoid by a finite set $S$. A
\emdef{weight} on $(G,S)$ is a function $\omega:S\to\R_+^*$. It induces
a \emdef{length} $\partial_\omega$ on $G$ by
$$\partial_\omega:\begin{cases}G\to\R_+\\
  g\mapsto\min\{\omega(s_1)+\dots+\omega(s_n)|\,s_1\cdots s_n=_G g\}.\end{cases}$$
A \emdef{minimal form} of $g\in G$ is a representation of $g$ as
a word of minimal length over $S$.  The \emdef{growth} of $G$ with
respect to $\omega$ is then
$$\gamma_\omega:\begin{cases}\R_+\to\R_+\\ n\mapsto \#\{g\in
G|\,\partial_\omega(g)\le n\}.\end{cases}$$

The function $\gamma:\R_+\to\R_+$ is \emdef{dominated} by
$\delta:\R_+\to\R_+$, written $\gamma\precsim\delta$, if there is a
constant $C\in\R_+$ such that $\gamma(n)\leq\delta(Cn)$ for all
$n\in\R_+$. Two functions $\gamma,\delta:\R_+\to\R_+$ are
\emdef{equivalent}, written $\gamma\sim\delta$, if
$\gamma\precsim\delta$ and $\delta\precsim\gamma$.

The following lemmata are well known:
\begin{lem}
  Let $S$ and $S'$ be two finite generating sets for the group $G$,
  and let $\omega$ and $\omega'$ be weights on $(G,S)$ and $(G,S')$
  respectively. Then $\gamma_\omega\sim\gamma_{\omega'}$.
\end{lem}
\begin{proof}
  Let $C=\max_{s\in S}\partial_{\omega'}(s)/\omega(s)$. Then
  $\partial_{\omega'}(g)\leq C\partial_\omega(g)$ for all $g\in G$,
  and thus $\gamma_\omega(n)\le \gamma_{\omega'}(Cn)$, from which
  $\gamma_\omega\precsim\gamma_{\omega'}$. The opposite relation holds
  by symmetry.
\end{proof}

The \emdef{growth type} of a finitely generated group $G$ is the
$\sim$-equivalence class containing its growth functions; it will be
denoted by $\gamma_G$.

Note that all exponential functions $b^n$ are equivalent, and
polynomial functions of different degree are inequivalent; the same
holds for the subexponential functions $e^{n^\alpha}$. We have
$$0\precnsim n\precnsim n^2\cdots\precnsim e^{\sqrt n}\precnsim\dots\precnsim e^n.$$
Note also that the ordering $\precsim$ is not linear.

\begin{lem}
  Let $G$ be a finitely generated group. Then $\gamma_G\precsim e^n$.
\end{lem}
\begin{proof}
  Choose for $G$ a finite generating set $S$, and define the weight
  $\omega$ by $\omega(s)=1$ for all $s\in S$. Then
  $\gamma_\omega(n)\le |S|^n$ for all $n$, so $\gamma_G\precsim e^n$.
\end{proof}

If there is a $d\in\N$ such that $\gamma_G\precsim n^d$, the group $G$
is of \emdef{polynomial growth} of degree at most $d$; if
$\gamma_G\sim e^n$, then $G$ is of \emdef{exponential growth};
otherwise $G$ is of \emdef{intermediate growth}. The existence of
groups of intermediate growth was first shown by
Grigorchuk~\cite{grigorchuk:growth}.

\section{The Grigorchuk $2$-group}
Let $\Sigma^*$ be the set of finite sequences over $\Sigma=\{0,1\}$.
For $x\in\Sigma$ set $\overline x=1-x$. Define recursively the
following length-preserving permutations of $\Sigma^*$:
\begin{alignat*}{3}
  &&a(x\sigma)&=\overline x\sigma;\\
  b(0\sigma)&=0a(\sigma),&&&b(1\sigma)&=1c(\sigma);\\
  c(0\sigma)&=0a(\sigma),&&&c(1\sigma)&=1d(\sigma);\\
  d(0\sigma)&=0\sigma,&&&d(1\sigma)&=1b(\sigma).
\end{alignat*}
Then $G$, the Grigorchuk
$2$-group~\cite{grigorchuk:burnside,grigorchuk:gdegree}, is the group
generated by $S=\{a,b,c,d\}$.  It is readily checked that these
generators are of order $2$ and that $V=\{1,b,c,d\}$ is a Klein group.

Let $H=V^G$ be the normal closure of $V$ in $G$. It is of index $2$ in
$G$ and preserves the first letter of sequences; i.e.\ $H\cdot
x\Sigma^*\subset x\Sigma^*$ for all $x\in\Sigma$. There is a map
$\psi:H\to G\times G$, written $g\mapsto (g_0,g_1)$, defined by
$0g_0(\sigma)=g(0\sigma)$ and $1g_1(\sigma)=g(1\sigma)$. As $H=\langle
b,c,d,b^a,c^a,d^a\rangle$, we can write $\psi$ explicitly as
$$\psi:\begin{cases}b\mapsto(a,c),\quad b^a\mapsto(c,a)\\ c\mapsto(a,d),\quad
  c^a\mapsto(d,a)\\ d\mapsto(1,b),\quad d^a\mapsto(b,1).\end{cases}$$

\section{The Growth of $G$}
Let $\eta\approx0.811$ be the real root of the polynomial
$X^3+X^2+X-2$, and define the following function on $S$:
\begin{alignat*}{2}
  \omega(a)&=1-\eta^3=\eta^2+\eta-1,&\qquad\omega(c)&=1-\eta^2, \\
  \omega(b)&=\eta^3=2-\eta-\eta^2,&\qquad\omega(d)&=1-\eta.
\end{alignat*}
It is a weight, because it takes positive values on every generator.

\begin{lem}\label{lem:nf}
  Every $g\in G$ admits a minimal form
  $$[*]a*a*a\dots*a[*],$$
  where $*\in\{b,c,d\}$ and the first and last $*$s are optional.
\end{lem}
\begin{proof}
  Clearly $\omega(s)>0$ for $s\in S$, so $\omega$ is a weight.  Let
  $w$ be a minimal form of $g$. The lemma asserts that one can suppose
  there are no consecutive letters in $\{b,c,d\}$ in $w$; now two
  equal consecutive letters cancel, and the product of any two
  distinct letters in $\{b,c,d\}$ equals the third one. For any
  arrangement $(x,y,z)$ of $\{b,c,d\}$ we have
  $\omega(x)\le\omega(y)+\omega(z)$, so the substitution of $z$ for
  $xy$ will not increase the weight of $w$.
\end{proof}

\begin{prop}\label{prop:shorten}
  Let $g\in H$, with $\psi(g)=(g_0,g_1)$. Then
  $$\eta\big(\partial_\omega(g)+\omega(a)\big)\ge\partial_\omega(g_0)+\partial_\omega(g_1).$$
\end{prop}
\begin{proof}
  Let $w$ be a minimal form of $g$. Thanks to
  Lemma~\ref{lem:nf} we may suppose the number of $*$s in $w$ is at
  most the number of $a$s plus one. Construct words $w_0,w_1$ over $S$
  using $\psi$ seen as a substitution on words; they represent $g_0$
  and $g_1$ respectively.  Note that
  \begin{align*}
    \eta(\omega(a)+\omega(b)) &= \omega(a)+\omega(c),\\
    \eta(\omega(a)+\omega(c)) &= \omega(a)+\omega(d),\\
    \eta(\omega(a)+\omega(d)) &= 0+\omega(b).
  \end{align*}
  As $\psi(b)=(a,c)$ and $\phi(aba)=(c,a)$, each $b$ in $w$
  contributes $\omega(a)+\omega(c)$ to the total weight of $w_0$ and
  $w_1$; the same argument applies to $c$ and $d$. Now, grouping
  together pairs of $*$s in $\{b,c,d\}$ and $a$s, we see that
  $\eta(\partial_\omega(g))$ is a sum of left-hand terms, possibly
  ${}-\eta\omega(a)$; while
  $\partial_\omega(g_0)+\partial_\omega(g_1)$ is bounded by the total
  weight of the letters in $w_0$ and $w_1$, which is precisely the sum
  of the corresponding right-hand terms.
\end{proof}

Let $\alpha=\log(2)/\log(2/\eta)\approx0.767$, and for $n\in\N$ set
$C_n=\frac1{n+1}\binom{2n}{n}$, the $n$th Catalan number; remember
that it is the number of labelled binary rooted trees with $n+1$
leaves~\cite[page~119]{catalan,vanlint-w:combinatorics}.
\begin{prop}\label{prop:subexp}
  Let $\zeta=\frac{\omega(a)}{2/\eta-1}$, let $K>\zeta$ be any
  constant, and for $n\in\R_+$ set
  $$L_n = \max\left\{1,\left\lceil2\left(\frac{n-\zeta}{K-\zeta}\right)^\alpha\right\rceil-1\right\}.$$
  Then we have
  \begin{equation}
    \label{eq:treebd}\gamma_\omega(n)\le C_{L_n-1}2^{L_n-1}\gamma_\omega(K)^{L_n}.
  \end{equation}
\end{prop}
\begin{proof}
  We construct an injection $\iota$ of $G$ into the set of labelled
  binary rooted trees each of whose leaves is labelled by an element
  of $G$ of weight bounded by $K$ and each of whose interior vertices
  is labelled by an element of the subgroup $\langle a\rangle$ of $G$.
  For $g\in G$, $\iota(g)$ is called its \emdef{representation}. It is
  constructed as follows: if $g\in G$ satisfies $\partial_\omega(g)\le
  K$, its representation is a tree with one vertex labelled by $g$. If
  $\partial_\omega(g)>K$, let $h\in\langle a\rangle$ be such that
  $gh\in H$, and write $\psi(gh)=(g_0,g_1)$. By
  Proposition~\ref{prop:shorten},
  $\partial_\omega(g_i)\le\eta\partial_\omega(g)$, so we may construct
  inductively the representations of $g_0$ and $g_1$. The
  representation of $g$ is a tree with $h$ at its root vertex and
  $\iota(g_0)$ and $\iota(g_1)$ attached to its two branches.
  
  We first claim that $\iota$ is injective: let $\tree$ be a tree in
  the image of $\iota$. If $\tree$ has one node labelled by $g$, then
  $\iota^{-1}(\tree)=\{g\}$. If $\tree$ has more than one vertex, let
  $h\in\langle a\rangle$ be the label of the root vertex and
  $(\tree_0,\tree_1)$ be the two subtrees connected to the root
  vertex. By induction on the number of vertices of $\tree$, we have
  $\tree_i=\iota(g_i)$ for unique $g_0$ and $g_1$.  Then as $\psi$ is
  injective there is a unique $g\in G$ with $\psi(gh)=(g_0,g_1)$, and
  $\iota^{-1}(\tree)=\{g\}$.

  We next prove by induction on $n$ that if $\partial_\omega(g)\le n$
  then its representation is a tree with at most $L_n$ leaves. Indeed
  if $n\le K$ then $g$'s representation has one leaf and $L_n=1$,
  while otherwise $g$'s representation is made up of those of $g_0$
  and $g_1$.  Say $\partial_\omega(g_0)=\ell$ and
  $\partial_\omega(g_1)=m$; then by Proposition~\ref{prop:shorten} we
  have $\ell+m\le\eta(n+\omega(a))$. By induction these
  representations have at most $L_\ell$ and $L_m$ leaves. As
  $\alpha<1$, we have $L_\ell+L_m\le2L_{(\ell+m)/2}$ for all $\ell,m$;
  and by direct computation, $L_{\eta/2(n+\omega(a))}=\lfloor
  L_n/2\rfloor$, so the number of leaves of $g$'s representation is
  $$L_\ell+L_m \le 2L_{(\ell+m)/2} \le 2L_{\eta/2(n+\omega(a))} \le L_n,$$
  as was claimed.
  
  We conclude that $\gamma(n)$ is bounded by the number of
  representations with $L_n$ leaves; there are $C_{L_n-1}$ binary
  trees with $L_n$ leaves, $2$ choices of labelling for each of the
  $L_n-1$ interior vertices, and $\gamma(K)$ choices for each leaf; so
  Equation~(\ref{eq:treebd}) follows.
\end{proof}

A lower bound on the growth of $G$ comes from the fact that $G$ is
residually a $2$-group:
\begin{thm}[Grigorchuk~\cite{grigorchuk:hp}]\label{thm:lowbd}
  Suppose $G$ is a finitely generated residually-$p$ group. Let
  $\{G_n\}$ be the Zassenhaus filtration of $G$. If $[G:G_n]\succnsim
  n^d$ for all $d$, then $\gamma_G\succsim e^{\sqrt n}$.
\end{thm}

\begin{proof}[Proof of Theorem~\ref{thm:main}]
  The sequence $[G:G_n]$ was shown to be of superpolynomial growth
  in~\cite{grigorchuk:hp}, so Theorem~\ref{thm:lowbd} yields the
  claimed lower bound; an elementary proof of this lower bound
  appears also in~\cite{grigorchuk:gdegree}.
  
  For the upper bound, which is the main result of the present note,
  we invoke Proposition~\ref{prop:subexp} with $K=1$, noting that
  $L_n\sim n^\alpha$ and $C_{L_n}\le4^{L_n}$, to obtain
  $\gamma_\omega\precsim (4\cdot2\cdot\gamma(1))^{n^\alpha}\sim
  e^{n^\alpha}$.
\end{proof}

\section{Conclusion}
The main fact used in the proof of Theorem~\ref{thm:main} is the
existence of minimal forms given by Lemma~\ref{lem:nf}, coming from
the natural map $\langle a\rangle*V\twoheadrightarrow G$. One can
impose stronger conditions on minimal forms, such as `not containing
$dada$ as a subword', coming from an explicit recursive presentation
of $G$~\cite{lysionok:pres}. Tighter upper bounds result from such
considerations. Yurij Leonov~\cite{leonov:lowerbd} recently obtained
improvements on the lower bound of Theorem~\ref{thm:main}.

I wish to thank Robyn Curtis, Igor Lys\"enok and Pierre de la Harpe
for having patiently listened to preliminary---and incorrect---proofs
of these results, and also of course Rostislav Grigorchuk without whom
these results couldn't even have existed.

\bibstyle{hamsalpha}
\bibliography{../mrabbrev,../people,../math,../bartholdi,../grigorchuk}

\def\nop#1{}\font\cyr=wncyr8
\providecommand{\bysame}{\leavevmode\hbox to3em{\hrulefill}\thinspace}
\begin{thebibliography}{Gro81}

\bibitem[Cat38]{catalan}
Eug{\`e}ne~C. Catalan, \emph{Note sur une \'equation aux diff\'erences finies},
  J. Math. Pures Appl. (9) \textbf{3} (1838), 508--516.

\bibitem[GH90]{ghys-h:gromov}
{\'E}tienne Ghys and Pierre~{de la} Harpe, \emph{Sur les groupes hyperboliques
  d'apr\`es {Mikhael Gromov}}, Progress in Mathematics, vol.~83, Birkh{\"a}user
  Boston Inc., Boston, MA, 1990, Papers from the Swiss Seminar on Hyperbolic
  Groups held in Bern, 1988.

\bibitem[Gri80]{grigorchuk:burnside}
Rostislav~I. Grigorchuk, \emph{On {B}urnside's problem on periodic groups},
  Funktsional. Anal. i Prilozhen. \textbf{14} (1980), no.~1, 53--54, English
  translation: {Functional Anal. Appl. \textbf{14} (1980), 41--43}.

\bibitem[Gri83]{grigorchuk:growth}
Rostislav~I. Grigorchuk, \emph{On the {M}ilnor problem of group growth}, Dokl.
  Akad. Nauk SSSR \textbf{271} (1983), no.~1, 30--33.

\bibitem[Gri84]{grigorchuk:gdegree}
Rostislav~I. Grigorchuk, \emph{Degrees of growth of finitely generated groups
  and the theory of invariant means}, Izv. Akad. Nauk SSSR Ser. Mat.
  \textbf{48} (1984), no.~5, 939--985, English translation: {Math. USSR-Izv.
  \textbf{25} (1985), no.~2, 259--300}.

\bibitem[Gri89]{grigorchuk:hp}
Rostislav~I. Grigorchuk, \emph{On the {H}ilbert-{P}oincar\'e series of graded
  algebras that are associated with groups}, Mat. Sb. \textbf{180} (1989),
  no.~2, 207--225, 304, English translation: {Math. USSR-Sb. \textbf{66}
  (1990), no.~1, 211--229}.

\bibitem[Gri91]{grigorchuk:kyoto}
Rostislav~I. Grigorchuk, \emph{On growth in group theory}, Proceedings of the
  International Congress of Mathematicians, Vol.\ I, II (Kyoto, 1990) (Tokyo),
  Math. Soc. Japan, 1991, pp.~325--338.

\bibitem[Gro81]{gromov:nilpotent}
Mikhael Gromov, \emph{Groups of polynomial growth and expanding maps}, Inst.
  Hautes {\'E}tudes Sci. Publ. Math. (1981), no.~53, 53--73.

\bibitem[Leo98]{leonov:lowerbd}
Yuri{\u\i}~G. Leonov, \emph{On lower estimation of growth for some torsion
  groups}, to appear, 1998.

\bibitem[LW92]{vanlint-w:combinatorics}
{Jacobus H. van} Lint and Richard~M. Wilson, \emph{A course in combinatorics},
  Cambridge University Press, 1992.

\bibitem[Lys85]{lysionok:pres}
Igor~G. Lysionok, \emph{A system of defining relations for the {Grigorchuk}
  group}, Mat. Zametki \textbf{38} (1985), 503--511.

\bibitem[Mil68]{milnor:solvable}
John~W. Milnor, \emph{Growth of finitely generated solvable groups}, J.
  Differential Geom. \textbf{2} (1968), 447--449.

\bibitem[Sva55]{svarts:growth}
A.~S. Svarts, \emph{A volume invariant of coverings}, Dokl. Akad. Nauk SSSR
  (1955), no.~105, 32--34 (Russian).

\bibitem[Tit72]{tits:linear}
Jacques Tits, \emph{Free subgroups in linear groups}, J. Algebra \textbf{20}
  (1972), 250--270.

\end{thebibliography}
\end{document}